# Design of Order-of-Addition Experiments

Jiayu Peng[1], Rahul Mukerjee[2] and Dennis K. J. Lin[1]

[1]Department of Statistics, The Pennsylvania State University, University Park, PA 16802, USA
[2]Indian Institute of Management Calcutta, Joka, Diamond Harbour Road, Kolkata 700 104, India

*Abstract*: In an order-of-addition experiment, each treatment is a permutation of $m$ components. It is often unaffordable to test all the $m$! treatments, and the design problem arises. We consider a model that incorporates the order of each pair of components and can also account for the distance between the two components in every such pair. Under this model, the optimality of the uniform design measure is established, via the approximate theory, for a broad range of criteria. Coupled with an eigen-analysis, this result serves as a benchmark that paves the way for assessing the efficiency and robustness of any exact design. The closed-form construction of a class of robust optimal fractional designs is then explored and illustrated.

## 1. Introduction

This paper considers a problem where the output of a process depends on the order of adding $m$ different components into the system, and interest is focused on understanding such dependence. Order-of-addition experiments have wide applications in, but not limited to, chemistry and related areas. For example, Ryberg (2008) studied a reaction in which the order of adding the reagents had a strong impact on the performance of the reaction, and a series of experiments was conducted to evaluate different orders of addition.

It is often unaffordable to test all the $m$! orders (for example, 10! is about 3.6 millions), and the design problem arises to choose a subset of orders for comparison. The naïve design which randomly chooses orders is often adopted in practice. For example, the properties of a phylogenetic tree depend on the order of addition of taxa (Olsen et al., 1994), and to our knowledge, often a number of randomly selected taxa orders are tested to reach conclusions (Stewart et al., 2001). Naturally, an optimal design, with a solid theoretical foundation, will be much more preferred than this naïve design. More order-of-addition related experiments can be found in Fuleki and Francis (1968), Shinohara and Ogawa (1998), and Ding et al. (2015), among others.

While the order-of-addition experiment is prevalent in practice and there is a strong need to investigate and improve upon its design, the statistical literature on this topic remains scarce. Van Nostrand (1995) appears to be the first statistical reference on the order-of-addition design. He suggested considering a design to detect the pairwise order effects that are often of much interest, with practitioners wishing to know, for example, if adding component 1 before 2 or vice versa has a significant influence on the response. Voelkel (2017) studied design criteria and design construction under a pairwise order model defined formally in § 2. A key feature of his work is the idea of ex-tending orthogonal arrays to designs that are naturally restricted, of which the order-of-addition orthogonal array is an important example. He showed that such an array leads to the same value of the *D*-criterion as the full design but did not report any optimality result.

We aim at initiating the systematic development of an optimal design theory for order-of-addition experiments. This is done under a model which is more flexible than the pairwise order model and covers it as a special case. Our results provide, in particular, a firm justification for Voelkel's (2017) work from the perspective of optimality and robustness. The tools that we employ include approximate theory via the use of signed permutation matrices, and an eigen-analysis motivated by an association algebra. All proofs appear in the appendix.

**2. Model formulation**

Suppose there are $m$ ($\geq 3$) components $1,\ldots,m$, which can be ordered in $m!$ ways. Any such ordering, say $a = a_1\ldots a_m$, which is a permutation of $1,\ldots,m$, is a treatment. As in Van Nostrand (1995) and Voelkel (2017), our model incorporates the order of every pair of components in any treatment. A new feature of the model is, however, that it allows for possible tapering, as may happen in practice, of the impact of any such pairwise order with an increase in the distance between the components in the pair. Thus, the impact of component $i$ preceding $j$ when they are at the two extremes of a treatment can be less severe than that, say, when they are next to each other.

To present the model formally, we write $A$ for the set of the $m!$ treatments. For $a \in A$, let $\tau(a)$ be the treatment mean of $a$, that is, the expectation of the response from $a$. As usual, it is assumed that the responses have equal variance and are uncorrelated. Let $S$ be the set of all pairs $ij$, $1 \leq i < j \leq m$. For $a$ ($=a_1\ldots a_m$)$\in A$ and $ij \in S$, write $h(ij,a)$ for the distance between $i$ and $j$ in $a$, that is, if $a_k = i$ and $a_l = j$, then $h(ij,a) = |k-l|$, so that $h(ij,a) \in \{1, \ldots, m-1\}$. Then according to our model

$$\tau(a) = \beta_0 + \Sigma_{ij} z_{ij}(a)\beta_{ij}, \qquad a \in A, \tag{1}$$

where $\Sigma_{ij}$ denotes sum over $ij \in S$, the $\beta_{ij}$ and $\beta_0$ are unknown parameters, and for each $ij \in S$,

$$z_{ij}(a) = c_{h(ij,a)}, \text{ if } i \text{ precedes } j \text{ in } a,$$

$$= -c_{h(ij,a)}, \text{ if } j \text{ precedes } i \text{ in } a. \tag{2}$$

Thus, if $m = 3$, then $\tau(312) = \beta_0 + c_1\beta_{12} - c_1\beta_{13} - c_2\beta_{23}$, and so on. In (2), the $c_h$, $h = 1,\ldots, m-1$, are known quantities such that $1 = c_1 \geq \ldots \geq c_{m-1} \geq 0$. These are intended to capture the possible tapering as indicated above. For example, one can take $c_h = 1/h$ or $c_h = c^{h-1}$, with known $c$, $0 < c < 1$. On the other hand, if no such tapering is anticipated then one can as well take $c_h = 1$ for all $h$, in which case (1) and (2) get reduced to the usual pairwise order model of Van Nostrand (1995) and Voelkel (2017). Our findings indicate a high degree of robustness of the optimal designs to the specific choice of the $c_h$. Hence one does not need to worry about such specifics at the design stage.

Let $\tilde{\beta} = (\beta_{12}, \beta_{13},\ldots, \beta_{m-1\,m})^T$, where the superscript denotes transpose. Then $\beta = (\beta_0, \tilde{\beta}^T)^T$ represents the parametric vector of interest. Similarly, let $z(a) = [z_{12}(a), z_{13}(a),\ldots, z_{m-1\,m}(a)]^T$ and $x(a) = [1, z(a)^T]^T$, $a \in A$. Then (1) can be expressed as

$$\tau(a) = \beta_0 + z(a)^T \tilde{\beta} = x(a)^T \beta. \qquad (3)$$

Write $q = m(m-1)/2$ and $p = q+1$. Then $\tilde{\beta}$ and $z(a)$ are $q \times 1$, while $\beta$ and $x(a)$ are $p \times 1$. Also,

$$x(a)^T x(a) = 1 + (m-1)c_1^2 + \cdots + c_{m-1}^2, \qquad (4)$$

for each $a$, because by (2), $m - h$ of the $z_{ij}(a)$, $ij \in S$, have absolute value $c_h$, $h = 1,\ldots, m-1$.

### 3. Optimality of the uniform design measure

Let $d_0$ be the full design which replicates each treatment once. Theorem 1, appearing later in this section, establishes the optimality of $d_0$ for inference on $\beta$ under a broad range of criteria, among all designs with the same number, $m!$, of runs. This result is rather strong because it holds irrespective of the $c_h$ in (2). Although $d_0$ is impractical for larger $m$, Theorem 1 will provide a very useful benchmark for assessing smaller designs.

One may anticipate the optimality of $d_0$ because, with run size $m!$, a design that replicates some treatments more than once while omitting some others altogether is intuitively unappealing. But, unlike in similar situations like traditional full factorials, the moment matrix of $d_0$, obtained later in § 4, is rather involved. For instance, it need not be completely symmetric (Kiefer, 1975) in the sense of having all diagonal elements equal as well as all off-diagonal elements equal. The $c_h$ in (2) further complicate matters. As a result, a direct combinatorial proof of the optimality of $d_0$ is quite challenging. An approach, based on the approximate theory is found to yield a subtle non-computational proof that does not even require explicit evaluation of the moment matrix of $d_0$.

To motivate the ideas, consider an $N$-run exact design $d$, where any treatment $a \in A$ is replicated $r(a)$ times and the integers $r(a)$ $(\geq 0)$ sum to $N$. By (3), $d$ has per run moment matrix

$$M(w) = \Sigma_a w(a)x(a)x(a)^T, \qquad (5)$$

where $w(a) = r(a)/N$ and $\Sigma_a$ denotes sum over $a \in A$. In approximate theory, the requirement on the design weights $w(a)$ to be integral multiples of $1/N$ is relaxed and these are allowed to be any nonnegative quantities, subject to $\Sigma_a w(a) = 1$. Then $w = \{w(a): a \in A\}$ is called a design measure, having moment matrix $M(w)$ as shown in (5). In particular, the full design $d_0$ corresponds to the uniform design measure $w_0$ over $A$ which has moment matrix

$$M_0 = M(w_0) = (1/m!)\Sigma_a x(a)x(a)^T. \qquad (6)$$

Let $\mathcal{M}$ denote the class of $p \times p$ nonnegative definite matrices. Recall that a signed permutation matrix is a square matrix having exactly one nonzero entry in each row and column, where any nonzero entry is 1 or –1. We consider optimality criteria $\phi(.)$ that are (i) concave over $\mathcal{M}$ and (ii) signed permutation invariant, that is, $\phi(R^T M R) = \phi(M)$ for every signed permutation matrix $R$ of order $p$ and every $M \in \mathcal{M}$. Given any such $\phi(.)$, a design measure is called $\phi$-optimal if it maximizes $\phi\{M(w)\}$ among all design measures. The commonly used $D$-, $A$- and $E$-criteria correspond to $\phi(M) = \log \det(M)$, $-\text{tr}(M^{-1})$ and $\lambda_{\min}(M)$, respectively, where $\log \det(M)$ and $-\text{tr}(M^{-1})$ are interpreted as $-\infty$ for singular $M$, and $\lambda_{\min}$ stands for the smallest eigenvalue. As all design measures here have the same $\text{tr}\{M(w)\}$ by (4) and (5), our framework also covers the $MS$-optimality criterion of Eccleston and Hedayat (1974), with $\phi(M) = -\text{tr}(M^2)$. This is equivalent to the $E(s^2)$ criterion for two-level supersaturated factorial designs (Booth and Cox, 1962) when $c_h = 1$ for all $h$, in which case each $x(a)$ has elements $\pm 1$. All these criteria and concave and signed permutation invariant.

Theorem 1 below establishes the optimality of the uniform design measure, or equivalently, that of the full design. Its proof exploits concavity and signed permutation invariance of the optimality criteria along with the fact, proved in the appendix, that the action of any permutation of $1,\ldots,m$ on the set of treatments induces an action of a signed permutation matrix on $\{x(a): a \in A\}$

**Theorem 1**. *The uniform design measure $w_0$ is $\phi$-optimal for every optimality criterion $\phi(.)$ which is concave and signed permutation invariant.*

Theorem 1 has several important implications as indicated below.

(i) It shows, in particular, the $D$-, $A$-, $E$- and $MS$-optimality of the uniform design measure $w_0$. By the equivalence theorem, the $D$-optimality of $w_0$ also implies $G$-optimality, that is,

$$\max\{x(a)^T M^{-1}(w) x(a): a \in A\} \geq p = \max\{x(a)^T M_0^{-1} x(a): a \in A\},$$

for every design measure $w$; see for example, Silvey (1980, Chapter 3). Thus, $w_0$ minimizes the maximum variance of the estimated responses at $a \in A$.

(ii) In view of (i), the full design $d_0$ is $D$-, $A$-, $E$-, $MS$- and $G$-optimal among all designs having $m!$ runs. An exact design, with a smaller number of runs but having the same moment matrix $M_0$ as $d_0$, also enjoys these optimality properties. More generally, such a design is $\phi$-optimal for every $\phi(.)$ as in Theorem 1. The matrix $M_0$ will be examined in more detail in § 4 with a view to assessing the efficiencies of any given exact design under the aforesaid criteria.

(iii) While $w_0$ may not be the unique design measure having the optimality properties mentioned in (i), any other design measure which is $D$-, $A$- or $MS$-optimal must have the same moment matrix $M_0$ as $w_0$. This is because these three criteria are strictly concave.

(iv) In particular, for the usual pairwise order model given by $c_h = 1$ for all $h$, order-of-addition orthogonal arrays of Voelkel (2017) have, by definition, moment matrix $M_0$, and therefore, as noted in (ii) above, these are, indeed, optimal in the sense of Theorem 1. Conversely, in this situation, an exact design that corresponds to a D-, A- or MS-optimal design measure must have moment matrix $M_0$, in view of (iii) above. Hence it is not hard to see that such a design has to be an order-of-addition orthogonal array. This settles an issue left open by Voelkel (2017).

## 4. Eigen-analysis and efficiency assessment

Although the proof of Theorem 1 does not require explicit knowledge of $M_0$, we need to find $M_0$ and its eigenvalues to assess the efficiencies of a given design measure or a given exact design under various optimality criteria. Some more notation will help. Write $I$ for the identity matrix of order $q$, and define $V$ as the $q \times q$ matrix, with rows and columns indexed by the elements of $S$, such that for $ij, kl \in S$, the $(ij,kl)$th element of $V$ is

$$V(ij,kl) = 1, \text{ if } i = k, j \neq l \text{ or } i \neq k, j = l,$$
$$= -1, \text{ if } i = l \text{ or } j = k, \qquad (7)$$
$$= 0, \text{ otherwise.}$$

For instance, if $m = 4$, then $S = \{12, 13, 14, 23, 24, 34\}$, and

$$V = \begin{bmatrix} 0 & 1 & 1 & -1 & -1 & 0 \\ 1 & 0 & 1 & 1 & 0 & -1 \\ 1 & 1 & 0 & 0 & 1 & 1 \\ -1 & 1 & 0 & 0 & 1 & -1 \\ -1 & 0 & 1 & 1 & 0 & 1 \\ 0 & -1 & 1 & -1 & 1 & 0 \end{bmatrix},$$

where $V(12,12) = V(12,34) = 0$, $V(13,23) = 1$, $V(23,12) = -1$, and so on. Let

$$b_0 = 2\{(m-1)c_1^2 + \cdots + c_{m-1}^2\}/\{m(m-1)\}, \qquad (8)$$

$$b_1 = 2\Sigma_h \{m - h(1) - h(2)\} c_{h(1)} \{2c_{h(1)+h(2)} - c_{h(2)}\}/\{m(m-1)(m-2)\}, \qquad (9)$$

where $\Sigma_h$ denotes sum over positive integers $h(1), h(2)$ such that $h(1) + h(2) \leq m - 1$. In general, (8) and (9) do not permit further simplification though, for $c_h = 1$ for all $h$, these reduce to

$$b_0 = 1, \quad b_1 = 1/3. \qquad (10)$$

**Theorem 2**. *The moment matrix of the uniform design measure $w_0$ is $M_0 = \text{diag}(1, b_0 I + b_1 V)$, having eigenvalues 1, $b_0 + (m-2)b_1$ and $b_0 - 2b_1$, with multiplicities 1, $m-1$ and $(m-1)(m-2)/2$, respectively.*

If the $-1$'s in $V$ were 1's, then $V$ would equal an association matrix of the triangular association scheme (Raghavarao, 1971, Chapter 8) and the related association algebra could be useful in studying the eigenvalues of $M_0$ in Theorem 2. Somewhat in the same spirit, the proof in the appendix obtains an expression for $V^2$ as a linear combination of $I$ and $V$, and makes use of it.

Theorem 2 enables us to assess the efficiencies of any design measure $w$, under various criteria, relative to the uniform design measure $w_0$ seen be optimal in Theorem 1. In particular, the $D$- and $A$-efficiencies of $w$, defined as $[\det\{M(w)\}/\det(M_0)]^{1/p}$ and $\text{tr}(M_0^{-1})/\text{tr}\{M^{-1}(w)\}$, are

$$D\text{-eff}(w) = \left[ \frac{\det\{M(w)\}}{\{b_0 + (m-2)b_1\}^{m-1}(b_0 - 2b_1)^{(m-1)(m-2)/2}} \right]^{1/p}, \quad (11)$$

and

$$A\text{-eff}(w) = \frac{1 + (m-1)[\{b_0 + (m-2)b_1\}^{-1} + \{(m-2)/2\}(b_0 - 2b_1)^{-1}]}{\text{tr}\{M^{-1}(w)\}}, \quad (12)$$

respectively. The design measure $w$ in (11) and (12) can well correspond to an exact design.

Equations (11) and (12) are pivotal in efficiency calculation for a given choice of the $c_h$ in (2) as well as studying robustness across such choices. For illustration, we study the performance of the order-of-addition orthogonal arrays of Voelkel (2017) under tapered models given by (i) $c_h = 1/h$, and (ii) $c_h = (1/2)^{h-1}$, $h = 1,\ldots, m-1$. Earlier, in § 3, these arrays were found to be optimal under the usual pairwise order model. By (11) and (12), for $m = 5$, the 12-run array in Voelkel's (2017) Table 3 has $D$- and $A$-efficiencies 0.985 and 0.972 under (i), and 0.989 and 0.980 under (ii). Similarly, for $m = 6$, the 24-run arrays in his Table 7 are also quite robust, the best of these being the leftmost one in that table, having $D$- and $A$-efficiencies 0.992 and 0.984 under (i), and 0.994 and 0.988 under (ii). Indeed, even these figures can be conservative as, under (i) or (ii), there may not exist any exact design of the same run size as these arrays and having moment matrix $M_0$

**5. Optimal fractional designs**

The current literature on order-of-addition designs lacks any systematic procedure for the construction of optimal fractional designs. To this end, we now propose a method that yields optimal fractions much smaller than the full design $d_0$, and has potential for further improvement as discussed later in § 6. Our construction is based on the fact, noted in § 3, that an $N$-run exact design $d$, having the same moment matrix $M_0$ as $d_0$, is optimal among all $N$-run designs, for every criterion $\phi(.)$ as in Theorem 1. Throughout this section, we work under the usual pairwise order model as given by $c_h = 1$ for all $h$. Then by (10) and Theorem 2,

$$M_0 = \text{diag}\{1, I + (1/3)V\}. \quad (13)$$

The resulting designs are seen to be very efficient under tapered models as well. We first illustrate the structure of the proposed designs, and then present the general construction procedure.

Let $m = 4$. Consider the 12-run half-fractional design $d$, as given in transposed form by

$$\begin{matrix} 1 & 2 & 4 & 3 & 1 & 3 & 4 & 2 & 1 & 4 & 3 & 2 \\ 2 & 1 & 3 & 4 & 3 & 1 & 2 & 4 & 4 & 1 & 2 & 3 \\ 3 & 4 & 1 & 2 & 2 & 4 & 1 & 3 & 2 & 3 & 1 & 4 \\ 4 & 3 & 2 & 1 & 4 & 2 & 3 & 1 & 3 & 2 & 4 & 1 \end{matrix},$$

This can be obtained from the first design in Table 2 of Voelkel (2017) by relabeling of components and row permutation. The design $d$ has moment matrix $M_0$, and hence enjoys the optimality properties of $d_0$. To see the structure of $d$, let

$$B_1 = \begin{bmatrix} 1 & 2 \\ 2 & 1 \end{bmatrix}, \quad \overline{B}_1 = \begin{bmatrix} 3 & 4 \\ 4 & 3 \end{bmatrix}, \quad B_2 = \begin{bmatrix} 1 & 3 \\ 3 & 1 \end{bmatrix}, \quad \overline{B}_2 = \begin{bmatrix} 2 & 4 \\ 4 & 2 \end{bmatrix}, \quad B_3 = \begin{bmatrix} 1 & 4 \\ 4 & 1 \end{bmatrix}, \quad \overline{B}_3 = \begin{bmatrix} 2 & 3 \\ 3 & 2 \end{bmatrix}.$$

Then $d$ can be expressed as $[D_1^T \ D_2^T \ D_3^T]^T$, where

$$D_u = \begin{bmatrix} B_u & \overline{B}_u \\ \sim \overline{B}_u & B_u \end{bmatrix}, \qquad u = 1, 2, 3. \tag{14}$$

with $\sim$ representing the operator of column reversal of any matrix; for example, $\sim \overline{B}_1$ has first column $(4, 3)^T$ and second column $(3, 4)^T$.

We now show that how the above design structure allows an extension to the case of general even $m$ ($\geq 4$). Let $s = m/2$, $L = m!/\{2(s!s!)\}$, and $\Gamma = \{1,\ldots, m\}$. Consider lexicographically arranged distinct sets $C_1,\ldots, C_L$, where each $C_u$ consists of 1 and $s - 1$ other elements of $\Gamma$. For each $u$, let $\overline{C}_u$ be the complement of $C_u$ in $\Gamma$. In both $C_u$ and $\overline{C}_u$, the elements are arranged in the ascending order. For $u = 1,\ldots, L$, let $B_u$ be the $s! \times s$ array with rows formed by all permutations of the elements of $C_u$, and define $\overline{B}_u$ similarly with reference to $\overline{C}_u$. Thus, if $m = 4$, then $C_1=\{1, 2\}$, $\overline{C}_1=\{3, 4\}$, $C_2=\{1, 3\}$, $\overline{C}_2=\{2, 4\}$, $C_3=\{1, 4\}$, $\overline{C}_3=\{2, 3\}$, entailing $B_u$ and $\overline{B}_u$, $u = 1, 2, 3$, as shown above. Similarly, if $m = 6$, then for example, $C_5 = \{1, 3, 4\}$ and $\overline{C}_5 = \{2, 5, 6\}$, and hence

$$B_5 = \begin{bmatrix} 1 & 1 & 3 & 3 & 4 & 4 \\ 3 & 4 & 1 & 4 & 1 & 3 \\ 4 & 3 & 4 & 1 & 3 & 1 \end{bmatrix}, \qquad \overline{B}_5 = \begin{bmatrix} 2 & 2 & 5 & 5 & 6 & 6 \\ 5 & 6 & 2 & 6 & 2 & 5 \\ 6 & 5 & 6 & 2 & 5 & 2 \end{bmatrix},$$

For $u = 1,\ldots, L$, now define $D_u$ as in (14). Along the lines of the 12-run design shown above for $m = 4$, let $d^*$ be design consisting of the $m!/s!$ treatments given by the rows of $D = [D_1^T \ldots D_L^T]^T$.

**Theorem 3**. *For every even $m \geq 4$, the design $d^*$ has moment matrix $M_0$ and is hence $\phi$-optimal, among designs with the same number of runs, for every optimality criterion $\phi(.)$ which is concave and signed permutation invariant.*

The above construction readily yields $\phi$-optimal designs for odd $m$ ($= 2s + 1$, $s \geq 2$) as well. To get such a design in $(2s +1)!/s!$ runs, one has to stack $2s +1$ copies of $D$. Then it suffices to insert a column consisting only of $2s+1$ just before the $l$th column of the $l$th copy ($l = 1,\ldots, 2s$), and also after the last column of the last copy.

In order to explore the performance of the above fractional designs under tapered models, we again consider (i) $c_h = 1/h$, and (ii) $c_h = (1/2)^{h-1}$, $h = 1,\ldots, m-1$. Quite reassuringly, the $D$- and $A$-efficiencies of our fractional designs for $m = 4,\ldots,10$, calculated from (11) and (12), turn out to be over 0.99, under both (i) and (ii). Because of the reason indicated at the end of § 4, even these figures can be conservative.

The optimal design $d^*$ in Theorem 3 permits a natural blocking, with treatments arising from each $D_u$ constituting one block. In the model (1) for the treatment mean $\tau(a)$, then $\beta_0$ has to be replaced by a block effect parameter depending on the block where any treatment $a$ in $d^*$ appears, and interest is focused on $\tilde{\beta} = (\beta_{12}, \beta_{13},\ldots, \beta_{m-1\,m})^{\mathrm{T}}$. The proof of Lemma A.4(a) in the appendix shows that the blocking of $d^*$ as envisaged above is a case of orthogonal blocking under the usual pairwise order model. Therefore, under this model, the resulting block design remains $\phi$-optimal for inference on $\tilde{\beta}$, for every monotone criterion $\phi(.)$ which is concave and signed permutation invariant, such as the $D$- and $A$-criteria. Moreover, for $m$ up to 10, the $D$- and $A$-efficiencies of these blocked fractional designs, calculated from appropriate versions of (11) and (12), are again found to be over 0.99, under both tapered models as given by (i) and (ii) of the last paragraph.

## 6. Discussion

Theorem 3 is a first step towards systematic construction of optimal fractional designs for order-of-addition experiments. There is a need to develop optimal or efficient designs in even smaller run sizes. Theorems 1 and 2, along with the resulting expressions (11) and (12) for $D$- and $A$-efficiencies which facilitate the study of robustness, are powerful tools for this purpose. Initial studies suggest that a refinement of the procedure in § 5, via the use of certain incomplete block designs and partial, rather than full, permutations, should work. Some progress has already been made in this direction covering, in particular, the cases of $m = 7,\ldots,12$. For example, with $m = 8$, we have found a design that requires only 168 runs and has moment matrix $M_0$ when $c_h = 1$ for all $h$. Thus, under the usual pairwise model, this design is optimal in the sense of Theorem 3. It is also very robust to the choice of the $c_h$, and has $D$- and $A$-efficiencies over 0.99 under tapered models as given by $c_h = 1/h$, and $c_h = (1/2)^{h-1}$, $h = 1,\ldots, m-1$. Step-down and exchange algorithms are quite promising in the construction of still smaller efficient designs. Thus, with $m = 8$, from the 168-run design mentioned above, we could obtain a 84-run design having $D$- and $A$-efficiencies 0.986 and 0.973 for $c_h = 1$; 0.973 and 0.946 for $c_h = 1/h$; and 0.972 and 0.945 for $c_h = (1/2)^{h-1}$, $h = 1,\ldots, m-1$. As noted in § 4, these efficiency figures are conservative. Moreover, further streamlining of these algorithms should yield even better results.

Our results do not require normality of the responses, but assume their homoscedasticity. This is violated if the response variance is treatment dependent. While theoretical results are then hard to obtain, algorithms in approximate theory can be readily employed to find a *D*- or *A*-optimal design measure, say $w^*$, at least for $m \leq 7$; see, for example, Torsney and Martin-Martin (2009). Numerical studies based on comparison with $w^*$ show that our optimal designs retain high efficiency when the variances do not differ widely. In contrast, if such heteroscedasticity is more severe then, as intuitively expected, $w^*$ turns out to be far from uniform, and this springs a pleasant surprise, making the conversion of $w^*$ to an efficient exact design rather easy. One needs to include only the treatments where $w^*$ assigns greater weights and exclude others. Such designs are seen to be quite robust, with high efficiency across various tapered models and over an appreciable range of the variances. Thus, with $m = 4$, if it is believed that responses from treatments beginning with component 1 may have a common variance greater than that of the rest, then the 12-run design {2134, 2143, 2341, 2431, 3124, 3142, 3241, 3421, 4123, 4132, 4231, 4321}, obtained as above, has *D*-efficiency over 0.98 and *A*-efficiency over 0.95, for variance ratio greater than 1.2 and all three choices of the $c_h$ as in the last paragraph. While a good design here is expected to exclude treatments beginning with 1, approximate theory guides us about which 12 of the remaining 18 should be included. Incidentally, these heuristics do not work under homoscedasticity where, by Theorem 1, the uniform design measure is optimal. This gives no immediate clue about which treatments should be retained in a smaller design and necessitates the development of further theory. Of course, the optimal design algorithms can become unmanageable for larger $m$, and more work is needed in this regard.

We conclude with the hope that the present endeavor will generate further interest in order-of-addition designs and related topics.

**Appendix**

**Proof of Theorem 1**

Let $\pi = \pi_1\ldots\pi_m$ be a permutation of $1,\ldots,m$, and let

$$\pi a = \pi_{a_1}\ldots\pi_{a_m}, \quad\quad\quad (A.1)$$

for any treatment $a = a_1\ldots a_m \in A$. Clearly, $A$ is also the set of all such $m!$ permutations $\pi$.

**Lemma A.1**. *Given any $\pi \in A$, there exists a signed permutation matrix $R(\pi)$ of order p such that $x(\pi a)^T = x(a)^T R(\pi)$, for every $a \in A$.*

*Proof.* Given any $\pi = \pi_1\ldots\pi_m$, for $ij \in S$, let $\bar{\pi}_i\bar{\pi}_j$ equal $\pi_i\pi_j$ if $\pi_i < \pi_j$, and $\pi_j\pi_i$ if $\pi_i > \pi_j$. By (A.1), for every $ij \in S$ and every $a \in A$, (i) $i$ precedes $j$ in $a$ if and only if $\pi_i$ precedes $\pi_j$ in $\pi a$, and (ii) the distance between $i$ and $j$ in $a$ equals that between $\pi_i$ and $\pi_j$ in $\pi a$. Recalling (2), therefore, by (i), $z_{\bar{\pi}_i \bar{\pi}_j}(\pi a)$ and $z_{ij}(a)$ have the same sign if and only if $\pi_i < \pi_j$, while by (ii), they have the same absolute value. In other words, $z_{\bar{\pi}_i \bar{\pi}_j}(\pi a)$ equals $z_{ij}(a)$ if $\pi_i < \pi_j$, and $-z_{ij}(a)$ if $\pi_i > \pi_j$. Thus, $z(\pi a)^T = z(a)^T \tilde{R}(\pi)$ for every $a \in A$, where $\tilde{R}(\pi)$ is a signed permutation matrix of order $q$ such that, for each $ij \in S$, the $(ij, \bar{\pi}_i \bar{\pi}_j)$th element of $\tilde{R}(\pi)$ is 1 or $-1$, according as $\pi_i < \pi_j$ or $\pi_i > \pi_j$, respectively, and all other elements of $\tilde{R}(\pi)$ are zeros. Hence $x(\pi a)^T = x(a)^T R(\pi)$, for every $a \in A$, where $R(\pi) = \text{diag}\{1, \tilde{R}(\pi)\}$ is a signed permutation matrix of order $p$. □

We now complete the proof of Theorem 1. Consider a design measure $w = \{w(a): a \in A\}$. For any $\pi \in A$, let $\pi w$ be the design measure that assigns, for each $a \in A$, weight $w(a)$ on treatment $\pi a$. Because $\phi(.)$ is concave, writing $\Sigma_\pi$ for sum over $\pi \in A$,

$$\phi\{(1/m!)\Sigma_\pi M(\pi w)\} \geq (1/m!)\Sigma_\pi \phi\{M(\pi w)\}. \tag{A.2}$$

Now, $\{\pi a: \pi \in A\} = A$, for every fixed $a \in A$, so that by (6), $\Sigma_\pi x(\pi a)x(\pi a)^T = \Sigma_a x(a)x(a)^T = m!M_0$. Hence by (5),

$$(1/m!) \Sigma_\pi M(\pi w) = (1/m!)\Sigma_\pi \Sigma_a w(a)x(\pi a)x(\pi a)^T = \Sigma_a w(a)M_0 = M_0. \tag{A.3}$$

Also, by (5) and Lemma A.1, any $\pi$ leads to a signed permutation matrix $R(\pi)$ such that

$$M(\pi w) = \Sigma_a w(a)x(\pi a)x(\pi a)^T = R(\pi)^T \{\Sigma_a w(a)x(a)x(a)^T\} R(\pi) = R(\pi)^T M(w) R(\pi),$$

and hence $\phi\{M(\pi w)\} = \phi\{M(w)\}$, because $\phi(.)$ is signed permutation invariant. Consequently, in view of (A.3), from (A.2) we get $\phi(M_0) \geq \phi\{M(w)\}$, and the result follows.

**Proof of Theorem 2**

**Lemma A.2**. (a) *For $ij \in S$, $\Sigma_a z_{ij}(a) = 0$ and $\Sigma_a z_{ij}^2(a) = (m!)b_0$.*

(b) *For $i, j, k$ satisfying $1 \leq i < j < k \leq m$,*

$\Sigma_a z_{ij}(a)z_{ik}(a) = (m!)b_1, \quad \Sigma_a z_{ik}(a)z_{jk}(a) = (m!)b_1, \quad \Sigma_a z_{ij}(a)z_{jk}(a) = -(m!)b_1.$

(c) *For $ij, kl \in S$, if the sets $\{i, j\}$ and $\{k, l\}$ are disjoint, then $\Sigma_a z_{ij}(a)z_{kl}(a) = 0$.*

*Proof.* (a) This follows from (2) and (8), noting that among $z_{ij}(a)$, $a \in A$, each of $c_h$ and $-c_h$ occurs with frequency $(m-2)!(m-h)$, $h = 1, \ldots, m-1$.

(b) By (2), considering the positions, in increasing order, occupied by some permutation of $i, j$ and $k$ in any treatment,

$\Sigma_a z_{ij}(a)z_{ik}(a) = 2\{(m-3)!\}\Sigma_u \{c_{u(2)-u(1)}c_{u(3)-u(1)} - c_{u(2)-u(1)}c_{u(3)-u(2)} + c_{u(3)-u(2)}c_{u(3)-u(1)}\},$

where $\Sigma_u$ denotes sum over integers $u(1), u(2), u(3)$ satisfying $1 \leq u(1) < u(2) < u(3) \leq m$. Given any positive integers $h(1), h(2)$ such that $h(1) + h(2) \leq m - 1$, there are $m - h(1) - h(2)$ triplets $\{u(1), u(2), u(3)\}$ as above such that $u(2) - u(1) = h(1)$ and $u(3) - u(2) = h(2)$. Therefore,

$$\Sigma_a z_{ij}(a) z_{ik}(a) = 2(m - 3)! \Sigma_h \{m - h(1) - h(2)\} \{c_{h(1)} c_{h(1)+h(2)} - c_{h(1)} c_{h(2)} + c_{h(2)} c_{h(1)+h(2)}\},$$

where $\Sigma_h$ is as in (9). Because

$$\Sigma_h \{m - h(1) - h(2)\} c_{h(2)} c_{h(1)+h(2)} = \Sigma_h \{m - h(1) - h(2)\} c_{h(1)} c_{h(1)+h(2)},$$

by symmetry, the first of the three identities in (b) is now evident from (9). The other two identities follow similarly.

(c) This is evident from (2), noting that the set of treatments $A$ can be partitioned into disjoint pairs such that the two treatments within each pair have the positions of $i$ and $j$ interchanged and every other component in the same position. □

**Lemma A.3.** (a) $V^2 = 2(m - 2)I + (m - 4)V$, (b) $V$ has eigenvalues $m - 2$ and $-2$, with respective multiplicities $m - 1$ and $(m - 1)(m - 2)/2$.

*Proof.* (a) For any $ij, kl \in S$, write $\rho$ for the $(ij, kl)$th element of $V^2$. It suffices to show that $\rho$ equals $2(m - 2)$ if $ij = kl$, and $(m - 4)V(ij, kl)$ otherwise. As $V$ is symmetric, $\rho$ equals the scalar product of the $ij$th and $kl$th rows of $V$. So, $\rho = \mu_{++} - \mu_{+-} - \mu_{-+} + \mu_{--}$, where $\mu_{+-}$ is the number of positions with the $ij$th row having 1 and the $kl$th row having $-1$, and so on. Part (a) now follows from (7) via consideration of the cases (i)-(vi) below. In what follows, $\mu = (\mu_{++}, \mu_{+-}, \mu_{-+}, \mu_{--})$.

  (i) If $ij = kl$, then $\mu = (m - i + j - 3, 0, 0, m + i - j - 1)$; $\rho = 2(m - 2)$.

  (ii) If $i = k$, $j \neq l$, then $i \leq m - 2$, $\mu = (m - i - 2, 0, 1, i - 1)$ or $(m - i - 2, 1, 0, i - 1)$, according as $j < l$ or $j > l$, respectively; $\rho = m - 4$.

  (iii) If $i \neq k$, $j = l$, then $j \geq 3$, $\mu = (j - 3, 1, 0, m - j)$ or $(j - 3, 0, 1, m - j)$, according as $i < k$ or $i > k$, respectively; $\rho = m - 4$.

  (iv) If $i = l$, then $i \geq 2$, $\mu = (1, m - i - 1, i - 2, 0)$; $\rho = -(m - 4)$.

  (v) If $j = k$, then $j \leq m - 1$, $\mu = (1, j - 2, m - j - 1, 0)$; $\rho = -(m - 4)$.

  (vi) If none of the five cases above arises, then $\{i, j\}$ and $\{k, l\}$ are disjoint sets. In this situation, $\mu$ equals $(1, 1, 1, 1)$ if $i > l$ or $k > j$, $(2, 2, 0, 0)$ if $i < k < l < j$, $(2, 0, 2, 0)$ if $k < i < j < l$, and $(2, 1, 1, 0)$ if $i < k < j < l$ or $k < i < l < j$. So, $\rho = 0$.

(b) By (a), any eigenvalue, $\lambda$, of $V$ satisfies $\lambda^2 = 2(m - 2) + (m - 4)\lambda$, and hence equals $m - 2$ or $-2$. As $tr(V) = 0$, these eigenvalues have multiplicities $m - 1$ and $(m - 1)(m - 2)/2$, respectively. □

We now complete the proof of Theorem 2. The expression for $M_0$ follows from (6), (7) and Lemma A.2, recalling that $x(a) = [1, z(a)^T]^T$. Hence by Lemma A.3, eigenvalues of $M_0$ are as stated.

**Proof of Theorem 3**

Let $X = [X_0, X_{12}, X_{13}, \ldots, X_{m-1\,m}]$ be the model matrix of $d^*$, arising from (1) and (2), with $c_h = 1$ for all $h$. Thus, $X_0$ is a column of ones and the $X_{ij}$ correspond to $\beta_{ij}$, $ij \in S$. Write $N = (m!)/(s!)$ for the run size of $d^*$. Because $d^*$ has moment matrix $(1/N)X^T X$, the result follows from (7), (13) and Lemma A.4 below, which is akin to Lemma A.2 but has a different proof.

**Lemma A.4**. (a) *For $ij \in S$, $X_0^T X_{ij} = 0$ and $X_{ij}^T X_{ij} = N$. Also, $X_0^T X_0 = N$.*

(b) *For $i, j, k$ satisfying $1 \le i < j < k \le m$, $X_{ij}^T X_{ik} = N/3$, $X_{ik}^T X_{jk} = N/3$ and $X_{ij}^T X_{jk} = -N/3$.*

(c) *For $ij, kl \in S$, if the sets $\{i, j\}$ and $\{k, l\}$ are disjoint, then $X_{ij}^T X_{kl} = 0$.*

*Proof.* (a) Clearly, $X_{ij}^T X_{ij} = X_0^T X_0 = N$, $ij \in S$, as $X$ has elements $\pm 1$. Next, let $G_u(ij)$ be the contribution of $D_u$ to $X_0^T X_{ij}$. By (2) and (14), $G_u(ij) = 0$, irrespective of whether $i$ is in $C_u$ or $\overline{C}_u$, and $j$ is in $C_u$ or $\overline{C}_u$. Hence $X_0^T X_{ij} = 0$, $ij \in S$.

(b) Write $G_u$ for the row vector with elements $G_u(ij, ik)$, $G_u(ik, jk)$ and $G_u(ij, jk)$, where $G_u(ij, ik)$ is the contribution of $D_u$ to $X_{ij}^T X_{ik}$, and $G_u(ik, jk)$, $G_u(ij, jk)$ are similarly defined. Because the rows of $B_u$ and $\overline{B}_u$ in (14) are formed by all permutations of the elements of $C_u$ and $\overline{C}_u$, respectively, in view of (7) and (13), the following hold.

   (i) $G_u = (2/3)(s!)(1, 1, -1)$, if either $i, j, k \in C_u$, or $i, j, k \in \overline{C}_u$;

   (ii) $G_u = 2(s!)(1, 0, 0)$, if either $j, k \in C_u$, $i \in \overline{C}_u$, or $i \in C_u$, $j, k \in \overline{C}_u$;

   (iii) $G_u = 2(s!)(0, 1, 0)$, if either $i, j \in C_u$, $k \in \overline{C}_u$, or $k \in C_u$, $i, j \in \overline{C}_u$;

   (iv) $G_u = 2(s!)(0, 0, -1)$, if either $i, k \in C_u$, $j \in \overline{C}_u$, or $j \in C_u$, $i, k \in \overline{C}_u$.

Because the situation (i) corresponds to $(2s-3)!/\{s!(s-3)!\}$ choices of $u$, and each of (ii), (iii), (iv) correspond to $(2s-3)!/\{(s-1)!(s-2)!\}$ choices of $u$, part (b) follows after a little algebra.

(c) Let $G_u(ij, kl)$ be the contribution of $D_u$ to $X_{ij}^T X_{kl}$. Then (2) and (14),

   (i) $G_u(ij, kl) = 2(s!)$, if either $i, k \in C_u$, $j, l \in \overline{C}_u$, or $j, l \in C_u$, $i, k \in \overline{C}_u$,

   (ii) $G_u(ij, kl) = -2(s!)$, if either $i, l \in C_u$, $j, k \in \overline{C}_u$, or $j, k \in C_u$, $i, l \in \overline{C}_u$,

and $G_u(ij, kl) = 0$ in all other situations. Part (c) is now immediate, because each of (i) and (ii) corresponds to $(2s-4)!/\{(s-2)!(s-2)!\}$ choices of $u$. □